\documentclass[final]{elsarticle}

\usepackage[ruled,vlined]{algorithm2e}
\usepackage{graphics}
\usepackage{amsmath,amssymb}
\usepackage{rotating}
\usepackage{multirow} 
\usepackage[usenames]{color}

\begin{document}

\begin{frontmatter}

\title{A Domain-decomposing parallel sparse linear system solver}

\author{Murat Manguoglu}

\address{Computer Engineering, Middle East Technical University, Ankara, Turkey, 06800 ({\tt manguoglu@ceng.metu.edu.tr}).}

\begin{abstract}
The solution of large sparse linear systems is often the most time-consuming part of many science and engineering applications. Computational fluid dynamics, circuit simulation, power network analysis, and material science are just a few examples of the application areas in which large sparse linear systems need to be solved effectively. In this paper we introduce a new parallel hybrid sparse linear system solver for distributed memory architectures that contains both direct and iterative components. We show that by using our solver one can alleviate the drawbacks of direct and iterative solvers,  achieving better scalability than with direct solvers and more robustness than with classical preconditioned iterative solvers. Comparisons to well-known direct and iterative solvers on a parallel architecture are provided.
\end{abstract}


\end{frontmatter}


\pagestyle{myheadings}
\thispagestyle{plain}
\markboth{M. MANGUOGLU}{A Domain Decomposing Parallel Sparse Linear System Solver}

\section{Introduction}

Many applications in science and engineering give rise to large sparse 
linear systems of equations. Some of these systems arise in the 
discretization of  partial differential equations (PDEs) modeling various 
physical phenomena, such as in computational fluid  dynamics ,semiconductor 
device simulations, and material science.  Large and sparse linear systems
also arise in applications that are not governed by PDEs (e.g. power  system 
networks, circuit simulation,  and graph problems ). 

\begin{figure}[tbp]
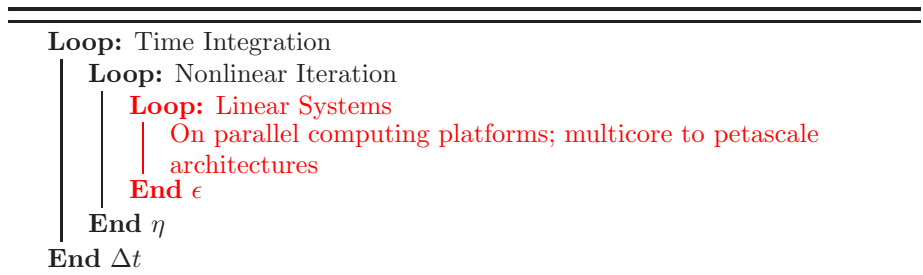
 
\begin{algorithm}[H]
\restylealgo{boxed}
\SetVline
\SetKwBlock{Loop}{Loop:}{End}
\Loop(Time Integration) {
\Loop(Nonlinear Iteration) {
\color{red} {\Loop(Linear Systems) {
On parallel computing platforms;
multicore to petascale architectures
} ($\epsilon$) }
 } ($\eta$)  
} ($\Delta t$)
\end{algorithm}
\label{computational_loop}
\caption{Target computational loop}
\end{figure}

Numerical simulation processes often consist of many layers of computational 
loops (e.g. see Figure~\ref{computational_loop}). It is a well known fact that 
the cost of the solution process is almost always governed by the solution of 
the linear systems especially for large-scale problems.

The emergence of multicore architectures and highly scalable platforms motivates
the development of novel algorithms and techniques that emphasize concurrency and 
are tolerant of deep memory hierarchies, as opposed to minimizing raw FLOP counts. 
While direct solvers are reliable, they are often  memory-intensive for large 
problems and offer limited scalability. Iterative solvers, on the other hand, 
are more efficient but, in the absence  of robust preconditioners, lack reliability.

In this paper we introduce a parallel sparse linear system solver that is hybrid. We note that we 
are using the term ``hybrid'' to emphasize that our solver is using both direct and iterative techniques. 
We advocate that using  our solver in hybrid mode one can  alleviate the  drawbacks of  
direct and iterative solvers, i.e. achieving more scalability  than a direct solver 
and more robustness than a classical preconditioned iterative solver.    

The rest of this paper is organized as follows. In Section 2, we discuss background and related work. 
In Section 3,  we give a description of the new algorithm and a simple example to demonstrate 
the details of the implementation. In Section 4, we present variety of numerical experiments. 
Finally, we conclude the paper with discussions in Section 5. 

\section{Background and related work}

Considerable effort has been spent on algebraic parallel sparse linear system solvers.  Sparse linear 
system solvers are traditionally divided into two groups (i) direct solvers (ii) iterative solvers. 
In the first group some examples are  MUMPS~\cite{1141208,587825,Amestoy00multifrontalparallel} , 
Pardiso~\cite{Schenk2004475,Schenk2006}, and SuperLU~\cite{superlu}.  

Iterative solvers mainly consist of classical preconditioned Krylov Subspace methods, 
and preconditioned Richardson iterations. Unlike direct sparse system solvers, iterative 
methods (with classical blackbox preconditioners) are not as robust. This is true even 
with the most recent advances  in creating LU-based preconditioners~\cite{benzi99orderings,
587288,ilupack}. Approximate inverse preconditioners~\cite{gravvanisfinite,gravvanis2009high,gravvanis2003solution,gravvanis2000explicit, benzi1996sparse}  are known to be more favorable for parallelism.  

The Spike algorithm~\cite{spike1,spike2,spike3,spike4,spike5,spike6,spike7} is a parallel solver 
for banded  systems,  that combines direct and iterative methods, is one of the first examples 
of hybrid  linear system solvers. More recently in~\cite{DBLP:journals/siamsc/ManguogluKSG10,manguoglu2,hyperthermia},
the Spike  algorithm was used for solving banded systems involving the preconditioner that is obtained 
after reordering the coefficient matrix with weights for sparse linear systems.

\section{Domain decomposing parallel solver}

We introduce a new parallel hybrid sparse linear system solver called Domain Decomposition 
Parallel Solver (DDPS) which can be used for solving sparse linear systems of equations: $Ax=f$. 
Recently,  we have presented an algorithm that used incomplete lu factorization for the 
diagonal block and its application on fluid structure interaction problems \cite{ijnmf, springerlink:10.1007/s00466-011-0619-0}. In this paper we 
introduce  DDPS that uses the direct solver Pardiso within each block and extend the results to general
 sparse systems from a variety of application areas.

We are motivated to create DDPS due to the  fact that many  applications use 
domain decomposition to distribute the work  among the  processors and the lack 
of reliability of black box preconditioned Krylov subspace methods and lack of scalability 
of direct solvers. METIS~\cite{Karypis98, Karypis99a} is often used  to partition the domain (and hence to partition the 
matrices).    DDPS is similar to  the Spike algorithm but unlike Spike it does not assume 
banded  structure for the  coefficient matrix $A$. Given a general sparse linear system  $Ax=f$,  we 
partition $A\in R^{n\times n}$ into $p$ block rows $A = [A_1 ,A_2, ..., A_p]^{T}$. Let
\begin{equation} 
A = \mathcal{D} + R, 
\end{equation}
where $\mathcal{D}$ consists of the $p$ block diagonals of $A$,
\begin{equation} 
\mathcal{D}= 
\begin{pmatrix}
A_{11}&        &        &   \\
      & A_{22} &        &   \\
      &        & \ddots &   \\
      &        &        &  A_{pp} \\
\end{pmatrix} 
\end{equation}
and $R$ consists of the remaining elements (i.e. $R = A - \mathcal{D}$). 
Let $\tilde{L}_i$ and $\tilde{U}_i$ be an incomplete LU factorizations 
of  $A_{ii}$ where $i=1,2,...,p$. We define 
\begin{equation}
\tilde{\mathcal{D}}=
\begin{pmatrix}
\tilde{A}_{11}&        &        &   \\
      & \tilde{A}_{22} &        &   \\
      &        & \ddots &   \\
     &        &        &  \tilde{A}_{pp} \\
\end{pmatrix}
\end{equation}
in which $\tilde{A}_{ii}=\tilde{L}_i\tilde{U}_i$. 

The DDPS algorithm is  shown in Figure~\ref{sspike}.   We assume the system $Ax=f$ is the one after METIS reordering 
\begin{figure}[htb]
\begin{algorithm}[H]
\restylealgo{boxed}
\SetVline
\KwData{$Ax=f$ and a partitioning information}
\KwResult{$x$} \SetKwBlock{Solve}{solve}{end}
1.\ $\mathcal{D}+R\longleftarrow A$ for the given  partitioning information\;
2.\ $\tilde{L}_i\tilde{U}_i \longleftarrow A_{ii}$ (approximate or exact)  for $i=1,2,...,p$ \;
3.\ $\tilde{R}\longleftarrow R$ (by dropping some elements) \;
4.\ $G\longleftarrow \tilde{\mathcal{D}}^{-1}\tilde{R}$ \; 
5.\ identify nonzero columns of $G$ and store their indices in array $c$ \; 
6.\ Solve $Ax=f$ via a Krylov subspace method with a preconditioner $P=\tilde{\mathcal{D}}+\tilde{R}$ and stopping tolerance $\epsilon_{out}$)
\end{algorithm}
\caption{DDPS algorithm.} 
\label{sspike} 
\end{figure}

\begin{figure}[htb]
\begin{algorithm}[H]
\restylealgo{boxed}
\SetVline
 \SetKwBlock{Solve}{solve}{end}
 \Solve($Pz=y$ ){
 ($\tilde{\mathcal{D}}^{-1}Pz = \tilde{\mathcal{D}}^{-1}y \Rightarrow (I + G)z =g$) \;
 \ 6.1 $g  \longleftarrow  \tilde{\mathcal{D}}^{-1}y$ \; 
 \ 6.2 $\hat{G} \longleftarrow  (I(c,c)+G(c,c))$; $\hat{z} \longleftarrow z(c)$; $\hat{g} \longleftarrow g(c)$ \;   
 \ 6.3 solve the smaller independent system: $\hat{G}\hat{z} = \hat{g}$ (directly or iteratively with stopping tolerance $\epsilon_{in}$)  \; 
 \ 6.4 $z(c) \longleftarrow \hat{z}$ \;
 \ 6.5 $z \longleftarrow g - Gz$  \;
 }
\end{algorithm}
\caption{Preconditioning operation: $Pz =y$ } 
\label{precond} 
\end{figure}

%
Stages 1-5 are considered as a preprocessing phase where the right hand side is not required. 
After preprocessing we solve the system via a Krylov subspace method and using a preconditioner. 
The major operations in a Krylov subspace method are: (i) matrix vector multiplications, (ii) inner 
products, and (iii) preconditioning operations in the form of $Pz=y$ (for some $y$). Only the details of the preconditioning 
operations for DDPS are given in Figure~\ref{precond}.

Each stage, with the exception of  solving the reduced system, can be executed with perfect parallelism requiring  
no interprocessor  communications.  The solution of the smaller system $\hat{G}\hat{z} = \hat{g}$ 
is the only  part of the preconditioning operation that require communication. The size of $\hat{G}$ is problem and partitioning  
dependent and it is expected to have an influence on the overall scalability of the algorithm.  The size of $\hat{G}$ is determined by the number of nonzero columns in $G$. We employ several 
techniques to reduce the size of $\hat{G}$: 
\begin{itemize} 
\item METIS reordering to reduce the total communication volume for given number of partitions and hence reducing the size of $\hat{G}$ by reducing the number of elements in $R$.  (We note that METIS works on undirected graphs, therefore we apply METIS on $(|A|+|A^T|)/2$).
\item A dropping strategy: Given a tolerance $\delta\in[0,1]$,  if for any column $k$ in 
$R_i$ $||R(:,k)_i||_\infty \leq \delta \times max_{j} ||R(:,j)_i||_\infty$ 
($i=1,2,..,p$) we do not consider that column when forming $\hat{G}$. 
Here $R_i$  is the block row partition of $R$ (i.e. $R = [R_1 , R_2, ..., R_p]^{T}$). Another possibility is to drop elements 
after computing $G$. In this paper, however,  we only consider the former as the latter is expected to be computationally expensive. 
\end{itemize} 

It is required the diagonal blocks, $A_{ii}$,  are nonsingular. In case they are singular, however,  in addition to METIS 
 applying HSL MC64 reordering and/or a diagonal perturbation can be considered. 

Notice that dropping elements from $R$ in stage 3 to reduce  the size of $\hat{G}$ results in 
an approximation of the solution. Furthermore,  we can use approximate LU-factorization 
of the diagonal blocks in stage 2 and solve $\hat{G}\hat{x} = \hat{g}$ iteratively in stage 6.3. 
Therefore, we place an  outer iterative layer  where we use the above algorithm 
as a solver for  systems involving the preconditioner  $P= \tilde{\mathcal{D}} + \tilde{R}$ 
where $\tilde{R}$  consists only of  the columns that are not dropped. We stop the outer iterations when 
the relative  residual at the $k^{th}$ iteration $||r_k||_\infty/||r_0||_\infty \leq \epsilon_{out}$. 

DDPS is a direct solver  if (i) nothing is dropped from $R$, (ii) exact
LU factorization of $A_{ii}$ is computed, and (iii) $\hat{G}\hat{z} = \hat{g}$ is solved exactly.
In the case of using DDPS as a direct solver, an outer iterative scheme may not be required 
but recommended.  In this paper we use the direct solver Pardiso for computing LU factorization of the diagonal blocks. 

The choices we make in  stages 2,3, and 6.3,  result in a solver that can be as robust as a 
direct solver or as scalable as  an iterative solver, or anything in between. 
Notice that the outer iterative layer also benefits from our partitioning strategy as 
METIS  reduces the total communication volume in parallel  sparse matrix 
vector multiplications. 

We note that $\hat{G}$ consists  of dense columns within each partition which we store as a 
two dimensional array in memory and as a result  matrix vector 
multiplications can be done via level 2 BLAS~\cite{blas2,blas3} (or level 3 in case of multiple 
right hand sides). 

In order to illustrate the steps of the basic DDPS algorithm (without any approximations) we provide the following system, $Ax=f$, with 9 unknowns, 
\begin{equation} 
\begin{pmatrix} 
{\color{red} \textbf{0.2}} & {\color{red} \textbf{1.0}} & {\color{red} \textbf{-1}} & 0 & 0.01 & 0 & 0 & 0 & -0.01 \\ 
{\color{red} \textbf{0.01}} & {\color{red}  \textbf{0.3}}& {\color{red} \textbf{0}} & 0 & 0 & 0 & 0 & 0 & 0 \\
{\color{red}  \textbf{-0.1}} & {\color{red}  \textbf{0}} & {\color{red} \textbf{0.4}}& 0 & 0.3 & 0 & 0 & 0 & 0 \\
0 & 0 & 0 & {\color{red} \textbf{0.3}} & {\color{red} \textbf{0.6}} & {\color{red} \textbf{2}} & 0 & 0 & 0 \\
0 & -0.2 & 0 & {\color{red}  \textbf{0}} & {\color{red} \textbf{0.4}} & {\color{red} \textbf{0}} & 0 & 0 & 1.1 \\
0 & 0 & 0 & {\color{red} \textbf{-0.2}} & {\color{red} \textbf{0.1}} & {\color{red}  \textbf{0.5}} & 0 & 0 & 0 \\
1.2 & 0 & 0 & 0 & 0 & 0 & {\color{red} \textbf{0.4}} & {\color{red}  \textbf{0.02}} & {\color{red} \textbf{3.0}} \\
0 & 0 & 0 & 0 & 0 & 0 & {\color{red}  \textbf{2.0}} & {\color{red} \textbf{0.5}} & {\color{red} \textbf{0}} \\
0 & 0 & 0 & 0 & 0 & 0 & {\color{red} \textbf{0}} & {\color{red} \textbf{0.1}} & {\color{red} \textbf{0.6}} 
\end{pmatrix} 
\begin{pmatrix} 
x_1 \\
x_2 \\
x_3 \\
x_4 \\
x_5 \\
x_6 \\ 
x_7 \\
x_8 \\
x_9 
\end{pmatrix}
= 
\begin{pmatrix} 
1 \\
1 \\
1 \\
1 \\
1 \\
1 \\
1 \\
1 \\
1 
\end{pmatrix}  
\end{equation}
Block diagonal matrix D is indicated in red color (or bold in black and white) for 3 partitions where each partition is of size 3. 
After  premultiplying both sides with $D^{-1}$ from left we obtain the modified system, $(I+G)x=g$  (we do not need to  form $D^{-1}$ explicitly to compute $D^{-1}R$)   
\begin{equation} 
\begin{pmatrix} 
{\color{blue} \textbf{1}} &{\color{blue}\textbf{0}}& 0& 0 & {\color{blue} \textbf{-9.12}} & 0 & 0 & 0 & {\color{blue} \textbf{0.12}} \\ 
{\color{blue}\textbf{0}}& {\color{blue} \textbf{1}} & 0& 0 & {\color{blue} \textbf{0.304}} & 0 & 0 & 0 & {\color{blue} \textbf{0.004}} \\
0&0& 1& 0 & -1.53 & 0 & 0 & 0 & 0.03 \\
0  & 0.0909 & 0 &1 &  0 &  0 & 0 & 0 & -0.5 \\
{\color{blue}\textbf{0}}  & {\color{blue} \textbf{-0.5}}  & 0 &  0 &  {\color{blue} \textbf{1}}  &  0 & 0 & 0 & {\color{blue} \textbf{2.75}} \\
0 & 0.1364 & 0 & 0 &  0 &  1 & 0 & 0 & -0.75 \\
0.5172 & 0 & 0 & 0 & 0 & 0 & 1 &  0 &  0 \\
-2.069 & 0 & 0 & 0 & 0 & 0 &  0 & 1 & 0 \\
{\color{blue} \textbf{0.3448}}  & {\color{blue}\textbf{0}} & 0 & 0 & {\color{blue}\textbf{0}} & 0 & 0 &  0 &  {\color{blue}\textbf{1}} 
\end{pmatrix} 
\begin{pmatrix} 
\mathbf{\color{blue}x_1} \\
\mathbf{\color{blue}x_2} \\
x_3 \\
x_4 \\
\mathbf{\color{blue}x_5} \\
x_6 \\ 
x_7 \\
x_8 \\
\mathbf{\color{blue}x_9}  
\end{pmatrix}
= 
\begin{pmatrix} 
{\color{blue}\textbf{-2}} \\
{\color{blue}\textbf{3.4}} \\
2 \\
-3.1818 \\
{\color{blue}\textbf{2.5}} \\
0.2273 \\
-1.3103 \\
7.2414 \\
{\color{blue}\textbf{0.4598}}  
\end{pmatrix}.  
\end{equation} 
We note that unknowns 1, 2, 5, and 9 form a smaller independent reduced system (indicated in blue color or bold in black and white) , 
\begin{equation} 
\begin{pmatrix} 
{\color{blue} \textbf{1}} & {\color{blue} \textbf{0}} & {\color{blue} \textbf{-9.12}} & {\color{blue}\textbf{0.12}} \\
{\color{blue}\textbf{0}} & {\color{blue}\textbf{1}} & {\color{blue}\textbf{0.304}} & {\color{blue}\textbf{-0.004}} \\
{\color{blue}\textbf{0}} & {\color{blue}\textbf{-0.5}} & {\color{blue}\textbf{1}} & {\color{blue}\textbf{2.75}} \\
{\color{blue}\textbf{0.3448}}  & {\color{blue}\textbf{0}}  & {\color{blue}\textbf{0}}  & {\color{blue}\textbf{1}}  
\end{pmatrix} 
\begin{pmatrix} 
\mathbf{\color{blue}x_1}\\
\mathbf{\color{blue}x_2}\\
\mathbf{\color{blue}x_5}\\
\mathbf{\color{blue}x_9} 
\end{pmatrix} 
=
\begin{pmatrix}
{\color{blue}\textbf{-2}}\\
{\color{blue}\textbf{3.4}}\\
{\color{blue}\textbf{2.5}}\\
{\color{blue}\textbf{0.4598}}  
\end{pmatrix} 
\end{equation} which has the solution $[x_1, x_2, x_5, x_9]^T =[-3.2389,3.4413,-0.1151,1.5766]^T$.   
Finally, we can retrieve the solution of the system via 
\begin{equation} 
\begin{pmatrix}
x_1\\
x_2\\
x_3\\
x_4\\
x_5\\
x_6\\
x_7\\
x_8\\
x_9 
\end{pmatrix}
= 
\begin{pmatrix} 
-2 \\
3.4 \\
2 \\
-3.1818 \\
2.5\\
0.2273 \\
-1.3103 \\
7.2414 \\
0.4598  
\end{pmatrix} 
-
\begin{pmatrix}
0 & 0          & -9.12 & 0.12 \\
0 & 0          & 0.304 & 0.004\\
0 & 0          & -1.53 & 0.03 \\
0 & 0.0909     & 0     & -0.5 \\
0 & -0.5       & 0     & 2.75 \\
0 & 0.1364     & 0     &-0.75\\
0.5172 &0      & 0     & 0 \\
-2.069 &0      & 0     & 0 \\
0.3448  & 0    & 0     & 0 
\end{pmatrix} 
\begin{pmatrix} 
x_1\\
x_2\\
x_5\\
x_9
\end{pmatrix} 
\end{equation} 
and obtain $x = [-3.2389, 3.4413,1.7766,-2.7063, -0.1151, 0.9405, 0.365, 0.5402, 1.5766]^T$.  

\section{Numerical experiments} 

The set of problems is obtained from the University of 
Florida Sparse Matrix Collection~\cite{davis97}. We choose the largest 
nonsymmetric matrix from each application domain. The list of the 
matrices and their properties are given in Table~\ref{tab:matrices}.  For 
each matrix we generate the corresponding right hand-side using a solution 
vector of all ones to ensure that $f\in span(A)$. 
\begin{sidewaystable}[h]
  \centering
  \caption{Linear systems from the University of Florida Sparse Matrix Collection, {\it  n} , {\it  nnz}, and  {\it dd} stands for  matrix size , number of nonzeros, and the degree of diagonal dominance, respectively}
    \begin{tabular}{l|llll}
    System  &{\it  n}     & {\it nnz} & {\it dd}   & {\it problem domain}  \\
    \hline    ATMOSMODL & 1,489,752 & 10,319,760 & 0 & computational fluid dynamics \\
    HVDC2 & 189,860 & 1,339,638 & 0 &power network \\
    LANGUAGE & 399,130 & 1,216,334 & $6.2\times 10^{-4}$  &directed weighted graph  \\
    OHNE2 & 181,343 & 6,869,939 & $1.4\times 10^{-11}$ &semiconductor device simulation \\
    RAJAT31 & 4,690,002 & 20,316,253 & 0 &circuit simulation \\
    THERMOMECH\_DK & 204,316 & 2,846,228 & 0.32  &thermal \\
    TMT\_UNSYM & 917,825 & 4,584,801 & 1 &electromagnetic \\
    TORSO3 & 259,156 & 4,429,042 & $9.9\times 10^{-2}$ &2D/3D problem \\
    XENON2 & 157,464 & 3,866,688 & $8.2\times 10^{-2}$ &material science \\
    \end{tabular}
  \label{tab:matrices}
\end{sidewaystable}
All numerical experiments are performed on an Intel 
Xeon (X5560@2.8GHz) cluster with Infiband interconnection 
and 16GB memory per node. The number of MPI processes is equal to 
the number of cores used and is also equal to the number of 
partitions  for DDPS. 

In the following numerical experiments,  we use a variation of preconditioned BiCGStab~\cite{templates} as the outer iterative solver. The smaller reduced system $\hat{G}\hat{z} = \hat{g}$ is also solved iteratively via BiCGStab without preconditioning. 
For the iterative solvers, the outer iterations are terminated 
when the number of iterations  reaches $1,000$ or the relative 
residual meets the stopping criterion ($||f-Ax||_\infty/||f||_\infty \leq 10^{-5}$). 
Failures of the solvers are indicated by F1 or F2 when the solver runs out 
of memory or the final relative residual is larger than $10^{-5}$, 
respectively. We limit the maximum number of iterations to 100 and 
the stopping tolerance to $\epsilon_{in}=10^{-4}$ for the inner 
iterations of   DDPS solver.    

We use ILUPACK with the following parameters,  reorderings: weighted matching and AMD~\cite{242944}, 
droptol: $10^{-1}$ , estimate for the condition numbers of the factors: $50$, and 
an elbow space of $10$ which are recommended by the user guide for general sparse 
linear systems.   ILUPACK uses  GMRES(30)  with a variation of incomplete LU factorization based preconditioner.  
MUMPS and Pardiso has been used with their default parameters and  using METIS reordering.  

In Table~\ref{tab:directcomp} we present the total solve time for MUMPS, Pardiso 
, DDPS($\delta=0.9$) and ILUPACK.  For 5 systems (indicated by blue) out of 9,  DDPS is faster 
than MUMPS (for 16 MPI processors). In addition, DDPS is more robust than ILUPACK and almost as robust as MUMPS direct solver, 
using 16 partitions DDPS fails  only in 2 cases while ILUPACK and  MUMPS  fails in 5 and 1 case, respectively. 
DDPS never runs out of memory while MUMPS runs out of memory for one of the problems unless more than 8 partitions are used.

\begin{sidewaystable}[htbp]
  \centering
  \caption{Total solve times (in seconds) for MUMPS, Pardiso, and DDPS, and  ILUPACK}
    \begin{tabular}{l|rrrrr|r|rrrr|r}
          & MUMPS &       &       &       &       & Pardiso & DDPS  &       &       &       & ILUPACK \\
    MPI Processes & 1     & 2     & 4     & 8     & 16    & 1     & 2     & 4     & 8     & 16    & 1 \\
    \hline
    ATMOSMODL & {\color{red}\textbf{F1}} & {\color{red}\textbf{F1}} & {\color{red}\textbf{F1}} & {\color{red}\textbf{F1}} & 171.3 & 1291.0 & 391.6 & 781.0 & 149.1 & 100.7 & 13.6 \\
    HVDC2 & 1.5   & 1.6   & 1.4   & 1.5   & 1.9   & 2.0   & 1.4   & 1.6   & 6.9   & {\color{red}\textbf{F2}} & {\color{red}\textbf{F2}} \\
    LANGUAGE & 504.6 & 273.9 & {\color{red}\textbf{F2}} & {\color{red}\textbf{F2}} & {\color{red}\textbf{F2}} & 1191.3 & 124.7 & 15.2  & 6.4   & 2.0   & 3.4 \\
    OHNE2 & 42.5  & 27.3  & 19.3  & 13.2  & 8.5   & 43.9  & 21.2  & 9.4   & {\color{red}\textbf{F2}} & {\color{red}\textbf{F2}} & {\color{red}\textbf{F2}} \\
    RAJAT31 & 78.3  & 67.6  & 59.3  & 54.1  & 53.7  & 57.9  & 258.7 & 150.5 & 106.2 & 45.1  & {\color{red}\textbf{F2}} \\
    THERMO & 2.9   & 2.3   & 2.1   & 2.1   & 2.8   & 3.0   & 2.0   & 10.5  & 20.5  & 11.2  & 6.8 \\
    TMT\_UNSYM & 14.5  & 12.1  & 10.3  & 9.8   & 9.7   & 10.8  & 170.7 & 140.2 & 99.0  & 77.8  & {\color{red}\textbf{F2}} \\
    TORSO3 & 40.2  & 26.3  & 18.2  & 12.4  & 9.6   & 49.4  & 20.6  & 10.0  & 4.0   & 2.1   & 2.2 \\
    XENON2 & 13.5  & 8.2   & 6.1   & 4.5   & 4.2   & 14.7  & 14.9  & 7.6   & 3.9   & 2.9   & {\color{red}\textbf{F2}} \\
    \end{tabular}%
  \label{tab:directcomp}%
\end{sidewaystable}%

The speedup with respect to Pardiso solver using a single core is given in Table~\ref{tab:si}.  We note that two problems achieve superlinear speed improvement due to cache effects. 

\begin{table}[htbp]
  \centering
  \caption{Speedup of DDPS  compared to Pardiso}
    \begin{tabular}{r|rrrrr}
          & Pardiso & DDPS  &       &       &  \\
    MPI Processes & 1     & 2     & 4     & 8     & 16 \\
    \hline
    ATMOSMODL & 1.0  & 3.3  & 1.7  & 8.7  & 12.8 \\
    HVDC2 & 1.0  & 1.4  & 1.2  & 0.3  & {\color{red}\textbf{F2}} \\
    LANGUAGE & 1.0  & 9.6  & 78.2 & 185.7 & 609.4 \\
    OHNE2 & 1.0  & 2.1  & 4.8  & {\color{red}\textbf{F2}} & {\color{red}\textbf{F2}} \\
    RAJAT31 & 1.0  & 0.2  & 0.4  & 0.6  & 1.3 \\
    THERMO & 1.0  & 1.5  & 0.3  & 0.2  & 0.3 \\
    TMT\_UNSYM & 1.0  & 0.1  & 0.1  & 0.1  & 0.1 \\
    TORSO3 & 1.0 & 2.4  & 4.9  & 12.2 & 23.6 \\
    XENON2 & 1.0  & 1.0  & 1.9  & 3.7  & 5.1 \\
    \end{tabular}%
  \label{tab:si}%
\end{table}%

In Table~\ref{tab:iterations}, the number of outer BiCGStab iterations for DDPS is provided as one increases the 
number of partitions. With the exception of two cases, namely {\em hvdc2} and {\em thermomech\_dk}, the number of iterations 
depends weakly (less than linear) on the number of partitions (or MPI processes).

\begin{table}[htbp]
  \centering
  \caption{Number of outer BiCGStab iterations for DDPS and ILUPACK}
    \begin{tabular}{r|rrrr|r}
          & DDPS  &       &       &       & ILUPACK \\
    MPI Processes & 2     & 4     & 8     & 16    & 1 \\
    \hline
    ATMOSMODL & 18    & 18    & 21.5  & 21.5  & 26 \\
    HVDC2 & 0.5   & 12.5  & 260   & {\color{red}\textbf{F2}} & {\color{red}\textbf{F2}} \\
    LANGUAGE & 5     & 7     & 6     & 6     & 4 \\
    OHNE2 & 0.5   & 0.5   & {\color{red}\textbf{F2}} &{\color{red} \textbf{F2}} & {\color{red}\textbf{F2}} \\
    RAJAT31 & 71.5  & 86.5  & 106.5 & 99    & {\color{red}\textbf{F2}} \\
    THERMO & 84.5  & 248.5 & 752.5 & 856   & 31 \\
    TMT\_UNSYM & 89.5  & 192   & 212   & 294   & {\color{red}\textbf{F2}} \\
    TORSO3 & 8.5   & 10    & 8.5   & 8.5   & 5 \\
    XENON2 & 53    & 63    & 67    & 90.5  & {\color{red}\textbf{F2}} \\
    \end{tabular}%
  \label{tab:iterations}%
\end{table}%

\begin{table}[htbp]
  \centering
  \caption{Average number of inner BiCGStab iterations for DDPS}
    \begin{tabular}{r|rrrr}
    MPI Processes & 2     & 4     & 8     & 16 \\
    \hline
    ATMOSMODL & 0.5   & 3.28  & 0.5   & 4.74 \\
    HVDC2 & 0.5   & 3.12  & 14.93 & {\color{red}\textbf{F2}} \\
    LANGUAGE & 0.5   & 0.5   & 0.92  & 1 \\
    OHNE2 & 3.5   & 3.5   & {\color{red}\textbf{F2}} & {\color{red}\textbf{F2}} \\
    RAJAT31 & 3.54  & 3.16  & 7.18  & 4.95 \\
    THERMO & 1     & 1     & 2.93  & 3.7 \\
    TMT\_UNSYM & 13.2  & 3.32  & 12.14 & 15.32 \\
    TORSO3 & 4.32  & 2.9   & 3.35  & 4.53 \\
    XENON2 & 1     & 1     & 1     & 4.7 \\
    \end{tabular}%
  \label{tab:avg}%
\end{table}%
The average number of inner BiCGStab iterations is given in Table~\ref{tab:avg}. Since we make sure the reduced system size is small via various techniques described earlier, number of iterations are relatively small for all systems with weak dependence on number of processes.

In Table~\ref{tab:rho} we show the effect of varying the drop tolerance, $\delta$, while the number 
of partitions is fixed at 16. A small $\delta$ results in a variation of DDPS that is  more like a direct 
solver. Although this causes the number of iterations to decrease, it also increases the memory requirement 
and the solver runs out of memory. For small $\delta$  memory problem appears in two cases, 
namely {\em rajat31} and {\em atmosmodl}. In 5 cases the number of outer iterations decrease as we decrease $\delta$. In the remaining 
two cases the DDPS solver failed even though $\delta$ is set to be a small number.    
\begin{table}[h]
  \centering
  \caption{Number of outer BiCGStab iterations for DDPS for 16 MPI processes}
    \begin{tabular}{l|llllll}
    $\delta$   & 0.99  & 0.9   & 0.6   & 0.3   & 0.1   & 1.0E-5 \\
    \hline
    ATMOSMODL & 23.5  & 21.5  & 19    & {\color{red}\textbf{F1}}    & {\color{red} \textbf{F1}}    & {\color{red} \textbf{F1}} \\
    HVDC2 & {\color{red}\textbf{F2}}    & {\color{red}\textbf{F2}}    & {\color{red}\textbf{F2}}    & {\color{red}\textbf{F2}}    & {\color{red}\textbf{ F2}}    & 8 \\
    LANGUAGE & 7     & 6     & 6     & 4     & 2.5   & 1 \\
    OHNE2 & {\color{red}\textbf{F2}}    & {\color{red}\textbf{F2}}    & {\color{red}\textbf{F2}}    & {\color{red}\textbf{F2}}    & {\color{red}\textbf{F2}}    & {\color{red}\textbf{F2}} \\
    RAJAT31 & 99    & 99    & 99    & {\color{red}\textbf{F1}}    & {\color{red}\textbf{F1}}    & {\color{red}\textbf{F1}} \\
    THERMOMECH\_DK & 645.5 & 856   & {\color{red}\textbf{F2}}    & {\color{red}\textbf{F2}}    & {\color{red}\textbf{F2}}    & 414 \\
    TMT\_UNSYM & 294   & 294   & 231   & {\color{red}\textbf{F2}}    & {\color{red}\textbf{F2}}    & {\color{red}\textbf{F2}} \\
    TORSO3 & 8.5   & 8.5   & 8.5   & 6.5   & 4     & 1 \\
    XENON2 & 99    & 90.5  & 105.5 & {\color{red}\textbf{F2}}    & {\color{red}\textbf{F2}}    & 1 \\
    \end{tabular}
  \label{tab:rho}
\end{table}

\section{Conclusion}
We have introduced a new hybrid sparse linear system solver called DDPS. We have shown that our new sparse linear system solver is often faster than  direct solvers and more robust than classical preconditioned Krylov subspace methods.  DDPS is flexible as it can be used in a variety of configurations.  Depending on the solver for the diagonal blocks a new variation of the algorithm will arise. The choice we make for solving the inner reduced system further increases the number of possibilities. Although we have used METIS to show the application of the algorithm on general sparse systems, DDPS algorithm is ideally suited for problems in which the matrices are already distributed via domain decomposition to minimize interprocessor communication.

\section*{Acknowledgments}
The author would like to thank Ahmed Sameh, Ananth Grama, David Kuck, Eric Cox, Faisal Saied,  Henry Gabb, Kenji Takizawa, and Tayfun Tezduyar for the numerous discussions and for their support.   This work has been partially supported by the European Community\rq{}s  Seventh  Framework  Programme  (FP7/2007-2013) under grant agreement no: RI-261557 and METU BAP-08-11-2011-128 grant

\bibliographystyle{elsarticle-num}
\bibliography{strings,main,MM-patch}   

\end{document}